\documentclass{amsart}%
\usepackage{amsmath}
\usepackage{amsfonts}
\usepackage{amssymb}
\usepackage{graphicx}%
\providecommand{\U}[1]{\protect\rule{.1in}{.1in}}
\newtheorem{theorem}{Theorem}
\theoremstyle{plain}
\newtheorem{acknowledgement}{Acknowledgement}

\newtheorem{corollary}{Corollary}

\newtheorem{lemma}{Lemma}

\newtheorem{remark}{Remark}

\numberwithin{equation}{section}
\begin{document}
\title[$q-$Gaussian distributions.]{$q-$Gaussian distributions. Simplifications and simulations}
\author{Pawe\l \ J. Szab\l owski }
\address{Department of Mathematics and Information Sciences\\
Warsaw University of Technology\\
pl. Politechniki 1\\
00-661 Warszawa, Poland}
\email{pszablowski@elka.pw.edu.pl; pawel.szablowski@gmail.com}
\date{September 2008}
\subjclass[2000]{62E17, 60E05; Secondary 68U20, 65C05}
\keywords{$q-Gaussian$ distribution, orthogonal polynomials, $q-Hermite$ polynomials.
rejection method.}

\begin{abstract}
We present some properties of measures ($q-$Gaussian) that orthogonalize the
set of $q-$Hermite polynomials. We also present an algorithm for simulating
i.i.d. sequences of random variables having $q-$Gaussian distribution.

\end{abstract}
\maketitle

\section{Introduction}

The paper is devoted to recollection of known and presentation of some new
properties of a distribution called $q-$Gaussian. We propose also a method of
simulation of i.i.d. sequences drown from it.

$q-$Gaussian is in fact a family of distributions indexed by a parameter
$q\in\lbrack-1,1]$. It is defined as follows.\newline For $q=-1,$ it is a
discrete $2$ point distribution, which assigns values $1/2$ to $-1$ and
$1$.\newline For $q\in\left(  -1,1\right)  ,$ it has density given by%

\[
f_{H}(x|q)=\frac{\sqrt{1-q}}{2\pi\sqrt{4-(1-q)x^{2}}}\prod_{k=0}^{\infty
}\left(  (1+q^{k})^{2}-(1-q)x^{2}q^{k}\right)  \prod_{k=0}^{\infty}%
(1-q^{k+1}),
\]
for $\left\vert x\right\vert \leq\frac{2}{\sqrt{1-q}}$. In particular
$f_{H}\left(  x|0\right)  \allowbreak=\allowbreak\frac{1}{2\pi}\sqrt{4-x^{2}%
},$ for $|x|\leq2.$ Hence it is Wigner distribution with radius $2$. \newline
For $q=1,$ $q$\emph{-Gaussian} distribution is the Normal distribution with
parameters $0$ and $1.$

Below, we present plots of $f_{H}\left(  x|-.4\right)  $ in \emph{blue,
}$f_{H}\left(  x|.1\right)  $ in \emph{orange,} $f_{H}\left(  x|.8\right)  $
in \emph{red }and\emph{\ }standard\emph{\ normal density }in \emph{black}%

\begin{center}
\includegraphics[
natheight=4.802300in,
natwidth=7.507400in,
height=2.0098in,
width=3.5898in
]%
{KK7MK900.wmf}%
\end{center}

This family of distributions was defined first in the paper of M. Bo\.{z}ejko,
B. K{\"{u}}mmerer and R. Speicher in 1997 in \cite{BKS97} in noncommutative
probability context. Later (\cite{bryc1} ) it appeared in quite classical
context namely as a stationary distribution $P_{H}$ of discrete time random
field $\mathbf{X=}\left\{  X_{n}\right\}  _{n\in\mathbb{Z}}$ defined by the
following relationships: $\mathbb{E}\left(  X_{i}\right)  \allowbreak
=\allowbreak0,$ $\mathbb{E}\left(  X_{i}^{2}\right)  \allowbreak
=\allowbreak1,$ $i\in\mathbb{Z},$
\begin{equation}
\exists a\in\mathbb{R};\forall n\in\mathbb{Z}:\mathbb{E}\left(  X_{n}%
|\mathcal{F}_{\neq n}\right)  =a\left(  X_{n-1}+X_{n+1}\right)  ,~a.s.
\label{bryc1}%
\end{equation}
and
\begin{gather}
\exists A,B,C\in\mathbb{R};\forall n\in\mathbb{Z}:\mathbb{E}\left(  X_{n}%
^{2}|\mathcal{F}_{\neq n}\right)  =\label{bryc2}\\
A\left(  X_{n-1}^{2}+X_{n+1}^{2}\right)  +BX_{n-1}X_{n+1}+C,~a.s.,
\end{gather}
where $\mathcal{F}_{\neq m}:=\sigma\left(  X_{k}:k\neq m\right)  $. It turns
out that parameters $a,$ $A,$ $B,$ $C$ are related to one another in such a
way that there are two parameters $q\geq-1$ and $0<|\rho|<1$ and all others
can be expressed through them:
\begin{align*}
a  &  =\frac{\rho}{1+\rho^{2}},A=\frac{\rho^{2}\left(  1-q\rho^{2}\right)
}{\left(  \rho^{2}+1\right)  \left(  1-q\rho^{4}\right)  },\\
B  &  =\frac{\rho^{2}\left(  1-\rho^{2}\right)  \left(  1+q\right)  }{\left(
\rho^{2}+1\right)  \left(  1-q\rho^{4}\right)  },C=\frac{\left(  1-\rho
^{2}\right)  ^{2}}{1-q\rho^{4}}.
\end{align*}
Then, one proves that%
\begin{equation}
\forall n\in\mathbb{Z},k,i\geq1:\mathbb{E}\left(  H_{k}\left(  X_{n}|q\right)
|\mathcal{F}_{\leq n-i}\right)  =\rho^{ki}H_{k}\left(  X_{n-i}|q\right)
,~a.s., \label{_*}%
\end{equation}
where $\mathcal{F}_{\leq m}:=\sigma\left(  X_{k}:k\leq m\right)  $, (similarly
one defines $\mathcal{F}_{\geq m}:=\sigma\left(  X_{k}:k\geq m\right)  $ ) and
$H_{k}\left(  x|q\right)  \allowbreak:\allowbreak k\geq-1$ are $q-$Hermite
polynomials defined below. It turns out that for $q>1$ the one-dimensional
distribution of the process $\mathbf{X}$ is not defined by moments. This case
is treated separately (e.g. in in \cite{Szab} ).

As mentioned earlier, here we will consider only the case $|q|\leq1.$ We will
preserve notation and denote family of $q-$Gaussian distributions by
$P_{H}\left(  q\right)  $ or simply $P_{H}.$

It turns out that there is quite large literature where this distribution
appears and is used to model different phenomena. See e.g. \cite{HLHM},
\cite{MBJW}, \cite{Ans04}, \cite{BB05}, \newline\cite{Bryc et al.}. Besides
random field defined above models notions that first appeared in
noncommutative context and hence establishes a link between noncommutative and
classical probability theories.

\begin{remark}
In the literature there exists another family of distributions under the same
name. It appears in the context of (Boltzmann-Gibbs)-statistical mechanics.
See e.g.\cite{TMNT2008} for applications and review.

Both families are indexed by basically one parameter $q\in\lbrack-1,1],$ and
for $q=1$ both include ordinary $N\left(  0,1\right)  $ distribution.
\end{remark}

In the sequel we will use the following traditional notation used in so called
'$q$-series theory' $\left[  0\right]  _{q}=0,$ $[n]_{q}=1+q+\ldots
+q^{n-1};n\geq1,\left[  0\right]  _{q}!\allowbreak=\allowbreak1,\left[
n\right]  _{q}!\allowbreak=\allowbreak%
{\displaystyle\prod\limits_{i=1}^{n}}
\left[  i\right]  _{q}$, $\left[
\begin{array}
[c]{c}%
n\\
k
\end{array}
\right]  _{q}\allowbreak=\allowbreak\frac{\left[  n\right]  _{q}!}{\left[
k\right]  _{q}!\left[  n-k\right]  _{q}!},$ for $n\geq0,$ $k=0,\ldots,n$ and
$0 $ otherwise, and $\left(  a|q\right)  _{n}=\prod_{k=0}^{n-1}\left(
1-aq^{k}\right)  ,$ for $n=1,2,\ldots,\infty$ (so called Pochhammer symbol).
Sometimes $\left(  a|q\right)  _{n}$ will be abbreviated to $\left(  a\right)
_{n}$ if it will not cause misunderstanding. Notice that $\left(  q\right)
_{n}\allowbreak=\allowbreak\left(  1-q\right)  ^{n}\left[  n\right]  _{q}!,$ $%
\genfrac{[}{]}{0pt}{}{n}{k}%
_{q}\allowbreak=\allowbreak\frac{\left(  q\right)  _{n}}{\left(  q\right)
_{k}\left(  q\right)  _{n-k}}$ and that $\left[  n\right]  _{q},$ $\left[
n\right]  _{q}!$ and $\left[
\begin{array}
[c]{c}%
n\\
k
\end{array}
\right]  _{q}$ tend to $n,$ $n!$ and $\binom{n}{k}$ (Newton's symbol)
respectively as $q->1.$

\begin{remark}
Introducing new variable $z$ defined by the relationship $1-\frac{(1-q)x^{2}%
}{2}\allowbreak=\allowbreak\cos2\pi z$ we can express $q-$Gaussian density
through Jacobi $\theta$ functions defined e.g. in \cite{Whittaker}. Namely we
have for $z\in\lbrack-1/2,1/2]$%
\[
f_{H}\left(  2\sin(\pi z)/\sqrt{1-q}|q\right)  =C_{q}\theta_{3}\left(
z|q\right)  \theta_{2}\left(  z|q\right)  ,
\]
with $C_{q}\allowbreak=\allowbreak\frac{\sqrt{1-q}\left(  q|q^{2}\right)
}{2\pi q^{1/4}\left(  q^{2}|q^{2}\right)  }$ where $\theta_{3}\left(
z|q\right)  $ and $\theta_{2}\left(  z|q\right)  $ are so called third and
second Jacobi Theta functions.
\end{remark}

Let us introduce family of polynomials (called $q-$Hermite) satisfying the
following three-term recurrence relationship
\begin{equation}
H_{n+1}(x|q)=xH_{n}(x|q)-[n]_{q}H_{n-1}(x|q), \label{qHer}%
\end{equation}
with $H_{-1}(x|q)=0,$ $H_{0}\left(  x|q\right)  =1.$ Notice that $H_{n}\left(
x|0\right)  \allowbreak=\allowbreak U_{n}\left(  x/2\right)  ,$ $n\geq-1,$
where $U_{n}\left(  x\right)  $ are Chebyshev polynomials of the second kind
defined by
\begin{equation}
U_{n}\left(  \cos\theta\right)  \allowbreak=\allowbreak\sin(\left(
n+1\right)  \theta)/\sin\theta, \label{Czeb2}%
\end{equation}
and that $H_{n}\left(  x|1\right)  \allowbreak=\allowbreak H_{n}\left(
x\right)  ,$ $n\geq-1$ where $H_{n}\left(  x\right)  $ are ('probabilist')
Hermite polynomials i.e. polynomials orthogonal with respect to Gaussian
$N\left(  0,1\right)  $ measure.

It turns out that $q$\emph{-Gaussian} is the distribution with respect to
which $q$\emph{-Hermite} polynomials are orthogonal. This fact can be easily
deduced from (\ref{_*}).

Thus in particular using the condition
\[
\forall n\geq1;q\in(-1,1):\int_{-2/\sqrt{1-q}}^{2/\sqrt{1-q}}H_{n}\left(
x|q\right)  f_{H}\left(  x|q\right)  dx=0,
\]
we can get all moments of $q$\emph{-Gaussian} distribution. Hence in
particular we have $\mathbb{E}X^{2n+1}\allowbreak=\allowbreak0,$ $n\geq0,$
$\mathbb{E}X^{2}\allowbreak=\allowbreak1,$ $\mathbb{E}X^{4}\allowbreak
=\allowbreak2+q,$ $\mathbb{E}X^{6}\allowbreak=\allowbreak5+6q+3q^{2}+q^{3} $
if only $X\sim P_{H}\left(  q\right)  .$

The aim of this paper is to make $q$\emph{-Gaussian} distribution more
friendly by presenting an alternative form of the density $f_{H}$ for
$q\in(-1,1)$, more easy to deal with (in particular we find the c.d.f. of
$P_{H}$), and suggest a method of simulation of i.i.d. sequences having
density $f_{H}.$

\section{Expansion of $f_{H}$}

In this section we will prove the following expansion theorem:

\begin{theorem}
\label{expansion1}For all $\left\vert q\right\vert <1,$ $\left\vert
x\right\vert \leq2/\sqrt{1-q}$ we have%
\begin{equation}
f_{H}(x|q)=\frac{\sqrt{1-q}}{2\pi}\sqrt{4-\left(  1-q\right)  x^{2}}\sum
_{k=1}^{\infty}(-1)^{k-1}q^{\binom{k}{2}}U_{2k-2}\left(  x\sqrt{1-q}/2\right)
, \label{ezpansion}%
\end{equation}
where $\binom{n}{k}\allowbreak=\allowbreak\frac{n!}{k!\left(  n-k\right)  !}.
$
\end{theorem}

As a corollary we get expression for the c.d.f. function of $P_{H}.$

\begin{corollary}
\label{d_function} The distribution function of $q-$Gaussian distribution is
given by
\begin{subequations}
\label{cdf}%
\begin{gather}
F_{H}\left(  y|q\right)  =\frac{1}{2}+\frac{1}{\pi}\arcsin\left(  y\sqrt
{1-q}/2\right)  \allowbreak\label{cdf1}\\
+\allowbreak\frac{1}{2\pi}\sqrt{4-\left(  1-q\right)  y^{2}}\sum_{k=1}%
^{\infty}\left(  -1\right)  ^{k-1}q^{\binom{k}{2}}\left(  1+q^{k}\right)
\frac{U_{2k-1}\left(  y\sqrt{1-q}/2\right)  }{2k}. \label{cdf2}%
\end{gather}

\end{subequations}
\end{corollary}

Identity (\ref{ezpansion}) can be a source of many interesting identities,
which may not be widely known outside the circle of researchers working in
special functions.

\begin{corollary}
For all $q\in(-1,1)$

$i)\forall~x\in\lbrack-\frac{2}{\sqrt{1-q}},\frac{2}{\sqrt{1-q}}]$%
\begin{equation}
\left(  q\right)  _{\infty}\prod_{k=1}^{\infty}\left(  \left(  1+q^{k}\right)
^{2}-\left(  1-q\right)  x^{2}q^{k}\right)  =\sum_{k=1}^{\infty}%
(-1)^{k-1}q^{\binom{k}{2}}U_{2k-2}\left(  x\sqrt{1-q}/2\right)  ,
\label{uproszcz}%
\end{equation}
where the polynomials $U_{k}$ are defined by (\ref{Czeb2}).

In particular, we have:

$ii)$%

\[
\left(  -q\right)  _{\infty}\left(  q^{2}|q^{2}\right)  _{\infty}=\sum
_{k=1}^{\infty}q^{\binom{k}{2}},
\]
a particular case of so called Jacobi's 'triple product identity'.

$iii)$%
\[
\left(  q^{3}|q^{3}\right)  _{\infty}=1+\sum_{k=1}^{\infty}(-1)^{3k}\left(
q^{\binom{3k}{2}}+q^{\binom{3k+1}{2}}\right)  ,
\]

$iv)$%
\begin{equation}
\left(  q\right)  _{\infty}^{3}=1+\sum_{k=2}^{\infty}(-1)^{k+1}\left(
2k-1\right)  q^{\binom{k}{2}}, \label{3_potega}%
\end{equation}

$v)$%
\begin{equation}
\forall n\in\mathbb{N}:%
\genfrac{[}{]}{0pt}{}{2n}{n}%
_{q}=\sum_{k=1}^{n}\left(  -1\right)  ^{k-1}\left(  1+q^{k}\right)
q^{\binom{k}{2}}%
\genfrac{[}{]}{0pt}{}{2n}{n-k}%
_{q}, \label{suma}%
\end{equation}

$vi)$ $f_{H}\left(  x|q\right)  $ is bimodal for $q\in(-1,q_{0}),$ where
$q_{0}(\cong-.107)$ is the largest real root of the equation $\sum
_{k=0}^{\infty}(2k+1)^{2}q^{k(k+1)/2}\allowbreak=\allowbreak0$.
\end{corollary}

\begin{lemma}
\label{aproximation}For all $q\in(-1,1)$ and $n\geq4$ we have\newline i)
\begin{align*}
&  \underset{\left\vert x\right\vert <2/\sqrt{1-q}}{\sup}\left\vert
f_{H}\left(  x|q\right)  -\frac{\sqrt{1-q}}{2\pi}\sqrt{4-\left(  1-q\right)
x^{2}}(\sum_{k=1}^{n-1}(-1)^{k-1}q^{\binom{k}{2}}U_{2k-2}\left(  x\sqrt
{1-q}/2\right)  \right\vert \\
&  \leq\frac{n\left\vert q\right\vert ^{(n-1)(n-2)/2}}{\pi\left(
1-q^{2}\right)  ^{2}}.
\end{align*}
\newline ii)%
\begin{align*}
&  \underset{\left\vert x\right\vert <2/\sqrt{1-q}}{\sup}|F_{H}\left(
y|q\right)  -\frac{1}{2}-\frac{1}{\pi}\arcsin\left(  y\sqrt{1-q}/2\right) \\
&  -\frac{1}{2\pi}\sqrt{4-\left(  1-q\right)  y^{2}}\sum_{k=1}^{n-1}\left(
-1\right)  ^{k-1}q^{\binom{k}{2}}\left(  1+q^{k}\right)  \frac{U_{2k-1}\left(
y\sqrt{1-q}/2\right)  }{2k}|\\
&  \leq\frac{\left\vert q\right\vert ^{n(n-1)/2}+\left\vert q\right\vert
^{n\left(  n+1\right)  /2}}{2\pi(1-\left\vert q\right\vert ^{n})}%
\end{align*}

\end{lemma}

\begin{remark}
Using the assertion of the above corollary one can approximate the density
$f_{H}$ as well as function $F_{H}\left(  y|q\right)  \allowbreak
-\allowbreak\frac{1}{2}\allowbreak-\allowbreak\frac{1}{\pi}\arcsin\left(
y\sqrt{1-q}/2\right)  $of $P_{H}$ by expressions of the type $\sqrt{4-\left(
1-q\right)  x^{2}}\times polynomial\ in\ x$ with great accuracy. This
expression is simple to analyze, simulate and calculate interesting
characteristics. Of course one should be aware that for small values of $n$
$\frac{\sqrt{1-q}}{2\pi}\sqrt{4-\left(  1-q\right)  x^{2}}\allowbreak
\times\allowbreak(\sum_{k=1}^{n}(-1)^{k-1}q^{\binom{k}{2}}U_{2k-2}\left(
x\sqrt{1-q}/2\right)  $ is not nonnegative for all $(1-q)x^{2}\leq4$! To give
a scent of how many $n$'s are needed to obtain the given accuracy we solved
numerically (using program Mathematica) the equation%
\[
\frac{n\left\vert q\right\vert ^{(n-1)(n-2)/2}}{\pi\left(  1-q^{2}\right)
^{2}}=\varepsilon,
\]
for several values $q$ and $\varepsilon.$ Let us denote by $N\left(
q,\varepsilon\right)  $ the solution of this equation. We have
\begin{tabular}
[c]{llllll}%
$\varepsilon=$%
$\backslash$%
$q=$ & $.1$ & $.4$ & $.7$ & $.9$ & $.99 $.\\
$.01$ & $3.59$ & $4.97$ & $7.71$ & $14.93$ & $56.73$\\
$.001$ & $4.04$ & $5.67$ & $8.73$ & $16.53$ & $60.86$\\
$.0001$ & $4.76$ & $6.26$ & $9.61$ & $17.97$ & $64.70$%
\end{tabular}
. \newline We also performed similar calculations for equation
\[
\frac{\left\vert q\right\vert ^{n(n-1)/2}+\left\vert q\right\vert ^{n\left(
n+1\right)  /2}}{2\pi(1-\left\vert q\right\vert ^{n})}=\varepsilon,
\]
obtaining:%
\begin{tabular}
[c]{llllll}%
$\varepsilon=$%
$\backslash$%
$q=$ & $.1$ & $.4$ & $.7$ & $.9$ & $.99 $.\\
$.01$ & $2.3$ & $3.3$ & $5.1$ & $9.5$ & $33$\\
$.001$ & $2.8$ & $4.2$ & $6.3$ & $11.5$ & $39$\\
$.0001$ & $3.2$ & $4.73$ & $7.3$ & $13.3$ & $44$%
\end{tabular}
.
\end{remark}

\section{Simulation}

There is an interesting problem of quick simulation of i.i.d. sequences drawn
from $q-Gaussian$ distribution, using few realizations of i.i.d. standard
uniform variates. One possibility is the \emph{rejection method }(see for
example \cite{Devroy}). It is not optimal in the sense that it uses least
realizations of independent, uniform on $[0,1]$ variates. But as one can see
below it works.

To apply this method one has to compare density of the generated variates with
another density that has the property of being 'easy generated' or another
words i.i.d. sequences of variables having this control density are easily
obtainable. In the case of density $f_{H}$ such natural candidate is
$\frac{\sqrt{(1-q)(4-(1-q)x^{2})}}{2\pi}.$ However this density is unimodal,
while the densities $f_{H}$ for $q$ below certain negative value are bimodal.
This would lead to inefficient simulation method requiring many trial
observations to be generated from $\frac{\sqrt{(1-q)(4-(1-q)x^{2})}}{2\pi}$ to
obtain one observation from $f_{H}$ for sufficiently small $q.$ That is why we
decided to take as 'easy' density the following one:%

\[
f_{E}\left(  x|q\right)  =\frac{\sqrt{(1-q)(4-(1-q)x^{2})}\prod_{j=1}%
^{3}((1+q^{j})^{2}-(1-q)q^{j}x^{2})}{2\pi\left[  9\right]  _{q}\left[
5\right]  _{q}},
\]
defined for $x\in\left(  \frac{-2}{\sqrt{1-q}},\frac{2}{\sqrt{1-q}}\right)  $.
However to be sure that this distribution can be used one has to prove the
following inequalities presented by the following Lemma.

\begin{lemma}
\label{symul}For $-1<q<1$ and $x\in\left[  \frac{-2}{\sqrt{1-q}},\frac
{2}{\sqrt{1-q}}\right]  $ we have:%
\[
\frac{f_{H}\left(  x|q\right)  }{f_{E}\left(  x|q\right)  }\leq M\left(
q\right)  ,
\]
where \newline$M\left(  q\right)  \allowbreak=\allowbreak\left\{
\begin{array}
[c]{cc}%
(1+q)(1-q^{3})(1-q^{5})(1-q^{9})\prod_{k=4}^{\infty}(1-q^{2k})(1+q^{k}) &
if~q\in(0,1)\\
(1+q)(1-q^{3})(1-q^{5})(1-q^{9})\prod_{k=4}^{\infty}\left(  1+\left\vert
q\right\vert ^{k}\right)  ^{2}\left(  1-q^{k}\right)  & if~q\in(-1,0)
\end{array}
\right.  .$
\end{lemma}

Function $M\left(  q\right)  $ has the following plot
\begin{center}
\includegraphics[
trim=0.000000in 0.000000in -0.001502in 0.005314in,
natheight=4.830800in,
natwidth=7.507400in,
height=4.1384cm,
width=6.4207cm
]%
{KK7MK901.wmf}%
\\
{\protect\small Fig. 1. }$M(q)(0\div5)${\protect\small \ versus }%
$q~(-1\div1).${\protect\small \ }%
\label{M(q)}%
\end{center}

Now following \cite{Devroy} we can simulate sequences of independent random
variables with $q-Gaussian$ distribution. If $q=\pm1$ then such simulation is trivial.

For $q\in(-1,1)$ we use Lemma \ref{symul} and program Mathematica. We
generated sequence of independent random variables from density $f_{E}$ by
inversion method (see \cite{Devroy}), since $f_{E}(x|q)$ can be integrated
leading to cumulative distribution function (c.d.f.):
\begin{align*}
\int_{-2/\sqrt{1-q}}^{x}f_{E}\left(  y|q\right)  dy  &  =D\left(  q\right)
x\sqrt{(1-q)(4-(1-q)x^{2})}Q_{2}\left(  x^{3},q\right) \\
&  +C\left(  q\right)  \arcsin(x\sqrt{1-q}/2)\overset{def}{=}\allowbreak
F_{E}\left(  x|q\right)  ,
\end{align*}
for $\left\vert x\right\vert \leq\,\frac{2}{\sqrt{1-q}},$ where $Q$ denotes
quadratic polynomial in $x^{3}$ with coefficients depending on $q,$ while the
constants $D$ and $C$ are known functions of $q.$ Recall that the inversion
method requires solving numerically the sequence of equations $F_{E}\left(
x|q\right)  \allowbreak=\allowbreak r_{i},$ where $r_{i}$ are observations
drown from standard uniform distribution.

Since the function $F_{E}\left(  x|q\right)  $ is strictly increasing on its
support and its derivative is known, there are no numerical problems in
solving this equation. Due to efficient procedure 'FindRoot' of Mathematica
solving this equation is quick.

Now let us recall how rejection method works in case $q\in(-1,1).$

Applying algorithm described in \cite{Devroy}, the rule to get one realization
of random variable having density $f_{H}$ is as follows.

\begin{enumerate}
\item we generate two variables: $X\sim f_{E}\left(  .\right)  $ and $Y\sim
U\left(  0,1\right)  $

\item set $T=M\left(  q\right)  f_{E}\left(  X|q\right)  /f_{H}\left(
X|q\right)  .$

\item If $YT>1$ then set $Z=X$ otherwise repeat (1) and (2).
\end{enumerate}

To see how this algorithm works, we present two simulation results performed
(consisting of $2000$ simulations) with $q=-.8$ (red dots) and $q=.8$ (green dots).%

\begin{center}
\includegraphics[
natheight=4.622400in,
natwidth=7.507400in,
height=1.8334in,
width=2.9715in
]%
{KK7MK902.wmf}%
\\
{\protect\small Fig. 2. Simulation of i.i.d. sequences from }%
${\protect\small P}_{H}${\protect\small .}%
\label{symulacje}%
\end{center}

Unfortunately this algorithm turns out to be very inefficient for $q$ close to
$-1$, more practically less than say $-.85.$ One can see this by examining
Figure 1. Values of $M\left(  q\right)  $ are very large then, showing that
one needs very large number of observations from density $f_{E}$ to obtain one
observation from $f_{H}.$ Thus there is still an open question to generate
efficiently observations from $f_{H}$ for values close to $-1.$

One might be inclined to use formula (\ref{cdf}) and inversion method applied
to its finite approximation ad again using procedure 'FindRoot'. Well we
applied this idea to simulate $2000$ observations from $P_{H}$ for $q=-.95.$
It worked giving the following results:
\begin{center}
\includegraphics[
natheight=4.427800in,
natwidth=7.507400in,
height=1.7365in,
width=2.9369in
]%
{KK7MK903.wmf}%
\\
{\protect\small Fig3. i.i.d. sequence from }$P_{H}${\protect\small \ by
inversion method for }$q=-.95.$%
\label{symulacje 2}%
\end{center}

We used procedure 'FindRoot' of Mathematica. It worked as one can see however
it lasted quite a time to get the result.

Besides, when we tried to get $2000$ observations from $P_{H}$ for $q=-.97,$
numerical errors seemed to play an important role as one can notice judging
from black dots that appeared between levels $0$ and $-.5$ on the picture
below.
\begin{center}
\includegraphics[
natheight=4.414000in,
natwidth=7.507400in,
height=1.7521in,
width=2.9715in
]%
{KK7MK904.wmf}%
\\
{\protect\small Fig. 4. i.i.d. sequence from }$P_{H}${\protect\small \ by
inversion method for }$q=-.97.$%
\end{center}

\section{Proofs}

\begin{proof}
[Proof of Theorem \ref{expansion1}]Let us denote $z=\frac{x\sqrt{1-q}}{2}.$
Hence $\left\vert z\right\vert <1.$ We have%
\[
f_{H}\left(  z\right)  =\frac{\left(  q\right)  _{\infty}\sqrt{1-q}}{4\pi
\sqrt{1-z^{2}}}\prod_{k=0}^{\infty}\left(  \left(  1+q^{k}\right)  ^{2}%
-4z^{2}q^{k}\right)  .
\]
Now let us notice that%
\begin{align*}
\left(  1+q^{k}\right)  ^{2}-4z^{2}q^{k}  &  =\left(  1+q^{k}-2zq^{k/2}%
\right)  \left(  1+q^{k}+2zq^{k/2}\right) \\
&  =\left(  \left(  \sqrt{1-z^{2}}+iz\right)  ^{2}+q^{k}\right)  \left(
\left(  \sqrt{1-z^{2}}-iz\right)  ^{2}+q^{k}\right)  .
\end{align*}
Now notice that since $\left\vert z\right\vert <1,$ we see that $\left\vert
\sqrt{1-z^{2}}+iz\right\vert =1.$ Thus we can write $\sqrt{1-z^{2}}%
+iz=\exp\left(  i\theta\right)  $ where $\theta=\arcsin z$ and also
$\sqrt{1-z^{2}}=\cos\theta.$ Hence we can write%
\[
\left(  1+q^{k}\right)  ^{2}-4z^{2}q^{k}=\left(  1+e^{2i\theta}q^{k}\right)
\left(  1+e^{-2i\theta}q^{k}\right)  ,
\]
and consequently%
\[
f_{H}\left(  z|q\right)  =\frac{\left(  q\right)  _{\infty}\sqrt{1-q}}%
{4\pi\cos\theta}\left(  -e^{2i\theta}\right)  _{\infty}\left(  1+e^{-2i\theta
}\right)  \left(  -qe^{-2i\theta}\right)  _{\infty}.
\]
We will now use so called 'triple product identity' (see \cite{AAR}, Theorem
10.4.1., p.497) that states in our setting, that%
\begin{align*}
\left(  q\right)  _{\infty}\left(  -e^{2i\theta}\right)  _{\infty}\left(
-qe^{-2i\theta}\right)  _{\infty}  &  =\sum_{k=-\infty}^{\infty}\left(
-1\right)  ^{k}q^{\binom{k}{2}}\left(  -e^{2i\theta}\right)  ^{k}%
=\sum_{k=-\infty}^{\infty}q^{k\left(  k-1\right)  /2}e^{2ik\theta}\\
&  =1+\sum_{k=1}^{\infty}q^{k(k-1)/2}e^{2ik\theta}+\sum_{j=1}^{\infty
}q^{j\left(  j+1\right)  /2}e^{-2ij\theta}\\
&  =1+\sum_{k=1}^{\infty}q^{\binom{k}{2}}e^{2ik\theta}+\sum_{k=2}^{\infty
}q^{\binom{k}{2}}e^{-2i\left(  k-1\right)  \theta}.
\end{align*}
Now notice that
\[
\frac{\left(  1+e^{-2i\theta}\right)  }{\cos\theta}=\frac{2\cos\theta
}{e^{i\theta}\cos\theta}=2e^{-i\theta}.
\]
Hence,%
\begin{align*}
f_{H}\left(  z|q\right)   &  =\frac{\sqrt{1-q}}{2\pi}e^{-i\theta}\left(
1+\sum_{k=1}^{\infty}q^{\binom{k}{2}}e^{2ik\theta}+\sum_{k=2}^{\infty
}q^{\binom{k}{2}}e^{-2i\left(  k-1\right)  \theta}\right)  =\\
&  =\frac{\sqrt{1-q}}{2\pi}e^{-i\theta}\times\\
&  \left(  1+e^{2\theta}+\sum_{k=2}^{\infty}q^{\binom{k}{2}}\left(
\cos2k\theta+\cos2\left(  k-1\right)  \theta\right)  +i\sum_{k=2}^{\infty
}q^{\binom{k}{2}}\left(  \sin2k\theta-\sin2\left(  k-1\right)  \theta\right)
\right) \\
&  =\frac{\sqrt{1-q}}{\pi}\cos\theta+\frac{\sqrt{1-q}}{\pi}\sum_{k=2}^{\infty
}q^{\binom{k}{2}}\cos\left(  2k-1\right)  \theta.
\end{align*}
To return to variable $z$ we have to recall definition of Chebyshev
polynomials. Namely, we have%
\begin{align*}
\cos\left(  \left(  2k-1\right)  \theta\right)   &  =\cos\left(  \left(
2k-1\right)  \arcsin z\right) \\
&  =\cos\left(  \left(  2k-1\right)  \left(  \frac{\pi}{2}-\arccos z\right)
\right) \\
&  =\left(  -1\right)  ^{k+1}\sin\left(  \left(  2k-1\right)  \arccos z\right)
\\
&  =\left(  -1\right)  ^{k+1}\sqrt{1-z^{2}}U_{2k-2}\left(  z\right)  ,
\end{align*}
where $U_{n}(z)$ is the Chebyshev polynomial of the second kind. More
precisely we have here:%
\begin{equation}
U_{n}(z)=\sin\left(  (n+1)\arccos z\right)  /\sqrt{1-z^{2}}. \label{czebyszev}%
\end{equation}
It is well known, that sequence $\left\{  U_{n}\right\}  $ satisfies
three-term recurrence equation%
\[
U_{n+1}(z)-2xU_{n}(z)+U_{n-1}(z)=0,
\]
with $U_{-1}(z)=0,$ $U_{0}\left(  z\right)  =1,$ and can be calculated
directly (see \cite{Bell}, Theorem 7.2, p. 188) as in (\ref{Czeb2}). Thus we
have shown that%
\[
f_{H}\left(  z\right)  =\frac{\sqrt{1-q}}{\pi}\sqrt{1-z^{2}}\left(  \sum
_{k=1}^{\infty}(-1)^{k-1}q^{\binom{k}{2}}U_{2k-2}\left(  z\right)  \right)  ,
\]
or equivalently%
\[
f_{H}\left(  x|q\right)  =\frac{\sqrt{1-q}}{2\pi}\sqrt{4-\left(  1-q\right)
x^{2}}\left(  \sum_{k=1}^{\infty}(-1)^{k-1}q^{\binom{k}{2}}U_{2k-2}\left(
x\sqrt{1-q}/2\right)  \right)  .
\]

\end{proof}

\begin{proof}
[Proof of Corollary \ref{d_function}]We have $F_{H}\left(  y\right)
\allowbreak=\allowbreak\int_{-2/\sqrt{1-q}}^{y}f_{H}\left(  x\right)
dx\allowbreak.$ Now we change variable to $z\allowbreak=\allowbreak
x\sqrt{1-q}/2$ and use (\ref{ezpansion}). Thus \newline$F_{H}\left(  y\right)
\allowbreak=\allowbreak\frac{2}{\pi}(\int_{-1}^{y\sqrt{1-q}/2}\sqrt{1-z^{2}%
}dz\allowbreak+\allowbreak\sum_{k=2}^{\infty}\left(  -1\right)  ^{k-1}%
q^{\binom{k}{2}}\int_{-1}^{y\sqrt{1-q}/2}\sqrt{1-z^{2}}U_{2k-2}\left(
z\right)  dz).$ We use now the following, easy to prove, formulae:
\newline$\int_{-1}^{y}\sqrt{1-z^{2}}U_{2n}\left(  z\right)  dz\allowbreak
=\allowbreak\sqrt{1-y^{2}}\left(  U_{2n+1}\left(  y\right)  /\left(
4n+4\right)  -U_{2n-1}\left(  y\right)  /\left(  4n\right)  \right)  $ and
$\int_{-1}^{y}\sqrt{1-x^{2}}dx=\allowbreak\frac{\pi}{2}\allowbreak
+\allowbreak\frac{1}{2}\arcsin y\allowbreak+\allowbreak\frac{1}{2}%
y\sqrt{1-y^{2}}\allowbreak$ and after rearranging terms get (\ref{cdf}).
\end{proof}

\begin{proof}
[Proof of Corollary \ref{tozsamosci}]Let $f_{U}\left(  x|q\right)
\allowbreak=\allowbreak\frac{\sqrt{1-q}}{2\pi}\sqrt{4-(1-q)x^{2}}$ for
$\left\vert x\right\vert \sqrt{1-q}\leq2.$ Assertion $i)$ is obtained directly
after noting that
\[
f_{H}\left(  x|q\right)  =f_{U}(x|q)\prod_{k=1}^{\infty}\left(  \left(
1+q^{k}\right)  ^{2}-\left(  1-q\right)  x^{2}q^{k}\right)  \left(
1-q^{k}\right)  .
\]
Following (\ref{ezpansion}), we get
\[
f_{H}\left(  x|q\right)  =f_{U}\left(  x|q\right)  \left(  \sum_{k=1}^{\infty
}(-1)^{k-1}q^{\binom{k}{2}}U_{2k-2}\left(  x\sqrt{1-q}/2\right)  \right)  .
\]
$ii)$ and $iii)$ are obtained by inserting $x=0$ and $x=1/\sqrt{1-q}$ in
(\ref{ezpansion}) and canceling out common factors. From (\ref{ezpansion}) it
follows also that values $U_{2n}\left(  0\right)  $ and $U_{2n}\left(
1/2\right)  $ will be needed. Keeping in mind (\ref{czebyszev}) we see that
$U_{2n}(0)\allowbreak=\allowbreak\cos\left(  n\pi\right)  \allowbreak
=\allowbreak\left(  -1\right)  ^{n}$ and
\[
U_{2n}\left(  1/2\right)  =\left\{
\begin{array}
[c]{ccc}%
1 & if & n=3m\\
0 & if & n=3m+1\\
-1 & if & n=3m+2
\end{array}
\right.  .
\]
On the other hand we see that $\left.  \prod_{k=1}^{\infty}\left(
(1+q^{k})^{2}-(1-q)x^{2}q^{k}\right)  \allowbreak\prod_{k=0}^{\infty
}(1-q^{k+1})\right\vert _{x=0}\allowbreak=\allowbreak%
{\displaystyle\prod\limits_{k=1}^{\infty}}
\left(  1-q^{2k}\right)  \left(  1+q^{k}\right)  $ and $\left.  \prod
_{k=1}^{\infty}\left(  (1+q^{k})^{2}-(1-q)x^{2}q^{k}\right)  \allowbreak
\prod_{k=0}^{\infty}(1-q^{k+1})\right\vert _{x=1/\sqrt{1-q}}\allowbreak
=\allowbreak\prod_{k=1}^{\infty}\left(  (1+q^{k})^{2}-q^{k}\right)
\allowbreak\prod_{k=1}^{\infty}(1-q^{k})\allowbreak=\allowbreak\prod
_{k=1}^{\infty}\left(  1+q^{k}+q^{2k}\right)  \allowbreak\prod_{k=1}^{\infty
}(1-q^{k})\allowbreak=\allowbreak\allowbreak\prod_{k=1}^{\infty}\frac
{1-q^{3k}}{1-q^{k}}\allowbreak\prod_{k=1}^{\infty}(1-q^{k})\allowbreak
=\allowbreak\prod_{k=1}^{\infty}\left(  1-q^{3k}\right)  .$

$iv)$ Putting $x=\pm\frac{2}{\sqrt{1-q}}$ in (\ref{uproszcz}) we get
\[
\prod_{k=1}^{\infty}\left(  1-q^{k}\right)  ^{3}=\sum_{k=1}^{\infty}%
(-1)^{k-1}q^{\binom{k}{2}}U_{2k-2}\left(  1\right)  .
\]
Now recall that $U_{2k}\left(  1\right)  =2k+1.$

$v)$ To see this notice that $q-$Hermite polynomials are orthogonal with
respect to the measure with density $f_{H}$. Thus we have%
\[
\forall m>0:\int_{-2/\sqrt{1-q}}^{-2/\sqrt{1-q}}H_{m}\left(  x|q\right)
f_{H}\left(  x,q\right)  dx=0.
\]
Using (\ref{ezpansion}) know that $\forall m>0$\newline$\allowbreak
\int_{-2/\sqrt{1-q}}^{-2/\sqrt{1-q}}H_{m}\left(  x|q\right)  \sqrt{4-\left(
1-q\right)  x^{2}}\left(  1+\sum_{k=2}^{\infty}(-1)^{k+1}q^{\binom{k}{2}%
}U_{2k-2}\left(  x\sqrt{1-q}/2\right)  \right)  dx\allowbreak=0.$ Observing
that function $f_{H}$ is symmetric and $q-$Hermite polynomials of odd order
are odd functions, we deduce that above mentioned identities are trivial for
odd $m.$ Thus, let us concentrate on even $m.$ Introducing new variable
$z\allowbreak=\allowbreak x\sqrt{1-q}/2,$ and multiplying both sides of this
identity by $\left(  1-q\right)  ^{m/2}$ we get%
\begin{equation}
\forall m>1:\int_{-1}^{1}h_{m}\left(  z|q\right)  \left(  1+\sum_{k=2}%
^{\infty}(-1)^{k+1}q^{\binom{k}{2}}U_{2k-2}\left(  z\right)  \right)
\sqrt{1-z^{2}}dz=0, \label{calka}%
\end{equation}
where $h_{m}\left(  z|q\right)  =\left(  1-q\right)  ^{m/2}H_{m}\left(
2z/\sqrt{1-q}|q\right)  .$ Polynomials $h_{m}$ are called continuous
$q-$Hermite polynomials. It can be easily verified (following (\ref{qHer}))
that they satisfy the following three-term recurrence equation%
\[
h_{n+1}\left(  t|q\right)  =2th_{n}\left(  t|q\right)  -(1-q^{n}%
)h_{n-1}\left(  t|q\right)  ,
\]
with $h_{-1}=0,$ $h_{0}\left(  z|q\right)  =1.$ Moreover, it is also known
that (see e.g. \cite{AAR}):%
\begin{equation}
h_{n}\left(  \cos\theta|q\right)  =\sum_{k=0}^{n}%
\genfrac{[}{]}{0pt}{}{n}{k}%
\cos\left(  \left(  n-2k\right)  \theta\right)  . \label{nowyHer}%
\end{equation}
Let us change once more variables in (\ref{calka}) and put $z=\cos\tau.$ Then,
for $\forall m>1:\allowbreak\int_{0}^{\pi}h_{m}\left(  \cos\tau|q\right)
\times\allowbreak\left(  1+\sum_{k=2}^{\infty}(-1)^{k+1}q^{\binom{k}{2}%
}U_{2k-2}\left(  \cos\tau\right)  \right)  \allowbreak\times\allowbreak
\sin\tau d\tau=0,$ or
\[
\forall m>1:\int_{0}^{\pi}h_{m}\left(  \cos\tau|q\right)  \left(  \sum
_{k=1}^{\infty}(-1)^{k-1}q^{\binom{k}{2}}\sin\left(  2k-1\right)  \tau\sin
\tau\right)  d\tau=0
\]
Keeping in mind that $2\sin\left(  2k-1\right)  \tau\sin\tau=\cos\left(
2k-2\right)  \tau-\cos\left(  2k\right)  \tau,$ we see that%
\[
\forall m>1:\int_{0}^{\pi}h_{2m}\left(  \cos\tau|q\right)  \left(  \sum
_{k=1}^{\infty}(-1)^{k-1}q^{\binom{k}{2}}\left(  \cos\left(  2k-2\right)
\tau-\cos\left(  2k\right)  \tau\right)  \right)  d\tau=0.
\]
Now keeping in mind that $\int_{0}^{\pi}\cos2m\tau d\tau=0$ for $m>0$ we see
that%
\[
\forall m>1:\int_{0}^{\pi}h_{2m}\left(  \cos\tau|q\right)  \left(
\sum_{k=m+1}^{\infty}(-1)^{k+1}q^{\binom{k}{2}}\left(  \cos\left(
2k-2\right)  \tau-\cos\left(  2k\right)  \tau\right)  \right)  d\tau=0.
\]
On the other hand taking into account (\ref{nowyHer}) we see that%
\[
\int_{0}^{\pi}h_{2m}\left(  \cos\tau|q\right)  \cos2k\tau d\tau=\pi%
\genfrac{[}{]}{0pt}{}{2m}{m-k}%
_{q}%
\]
for $k=0,1,\ldots,m,.$ Hence we have (\ref{suma}).

$vi)$ Keeping in mind that $f_{H}$ is symmetric with respect to $x$ we deduce
that the point of change of modality of $f_{H}$ must be characterized by the
condition $f_{H}^{\prime\prime}\left(  0|q_{0}\right)  \allowbreak
=\allowbreak0.$ Calculating second derivative of the right hand side of
(\ref{ezpansion}) and remembering that $\left.  \left(  \sqrt{4-\left(
1-q\right)  x^{2}}\right)  ^{\prime}\right\vert _{x=0}\allowbreak
=\allowbreak0$ we end up with and equation $0\allowbreak=\allowbreak
-\frac{1-q}{2}\sum_{k=0}^{\infty}q^{k(k+1)/2}\allowbreak+\allowbreak
2(1-q)\allowbreak\times\allowbreak\sum_{k=1}^{\infty}\left(  -1\right)
^{k}k(k+1)q^{k(k+1)/2}(-1)^{k+1}.$ Now since $4k(k+1)+1\allowbreak
=\allowbreak\left(  2k+1\right)  ^{2}$ we get equation in $vi)$ defining
$q_{0}.$
\end{proof}

To prove Lemma \ref{aproximation} we need the following lemma.

\begin{lemma}
\label{numerical}Suppose $0<r<1$ and $n\geq3.$ Then
\[
\sum_{k\geq n}(2k-1)r^{\binom{k}{2}}\leq\frac{2nr^{n\left(  n-1\right)  /2}%
}{(1-r^{2})^{2}}.
\]

\end{lemma}

\begin{proof}
Recall that for $\left\vert \rho\right\vert <1$ we have: $\sum_{i\geq1}%
i\rho^{i-1}\allowbreak=\allowbreak\left(  \frac{1}{1-\rho}\right)  ^{\prime
}\allowbreak=\allowbreak\frac{1}{\left(  1-\rho\right)  ^{2}}$ and that
$\sum_{i=0}^{m}\rho^{i}=\frac{1-\rho^{m+1}}{1-\rho}.$ Thus we have
$\sum_{k\geq n}(2k-1)r^{\binom{k}{2}}\allowbreak\allowbreak=$\newline%
$\allowbreak r^{\left(  n-1\right)  \left(  n-4\right)  /2}\sum_{k\geq
n}\left(  2k-1\right)  r^{k\left(  k-1\right)  /2-\left(  n-1\right)  \left(
n-4\right)  /2}\allowbreak=\allowbreak r^{\left(  n-1\right)  \left(
n-4\right)  /2}\times$\newline$\sum_{k\geq n}\left(  2k-1\right)
r^{2k-2}r^{k\left(  k-1\right)  /2-\left(  n-1\right)  \left(  n-4\right)
/2-2(k-1)}.$ Now notice that $k\left(  k-1\right)  /2-\left(  n-1\right)
\left(  n-4\right)  /2-2(k-1)=\allowbreak\frac{1}{2}\left(  k-n\right)
\left(  k+n-5\right)  \allowbreak\geq\allowbreak0$ for $k\geq n\geq3.$ Hence
$\sum_{k\geq n}(2k-1)r^{\binom{k}{2}}\allowbreak\leq\allowbreak r^{\left(
n-1\right)  \left(  n-4\right)  /2}\sum_{k\geq n}\left(  2k-1\right)
r^{2k-2}\allowbreak\leq\allowbreak2r^{\left(  n-1\right)  \left(  n-4\right)
/2}\sum_{k\geq n}k\left(  r^{2}\right)  ^{k-1}\allowbreak=\allowbreak
2r^{\left(  n-1\right)  \left(  n-4\right)  /2}\frac{d}{dr^{2}}\left(
\frac{1}{1-r^{2}}-\frac{1-r^{2n}}{1-r^{2}}\right)  \allowbreak=\allowbreak
2r^{\left(  n-1\right)  \left(  n-4\right)  /2}\frac{d}{dr^{2}}\left(
\frac{\left(  r^{2}\right)  ^{n}}{1-r^{2}}\right)  \allowbreak=$%
\newline$\allowbreak2r^{\left(  n-1\right)  \left(  n-4\right)  /2}%
\frac{n\left(  r^{2}\right)  ^{n-1}-\left(  n-1\right)  \left(  r^{2}\right)
^{n}}{\left(  1-r^{2}\right)  ^{2}}\allowbreak\allowbreak\leq\allowbreak
\allowbreak2\frac{nr^{n\left(  n-1\right)  /2}}{\left(  1-r^{2}\right)  ^{2}%
}.$
\end{proof}

\begin{proof}
[Proof of Lemma \ref{aproximation}]$i)$

We have
\[
U_{n}(\cos\theta)=\frac{\sin\left(  \left(  n+1\right)  \theta\right)  }%
{\sin\theta}.
\]
From this fact we deduce that $\underset{\left\vert x\right\vert \leq1}{\sup
}\left\vert U_{n}\left(  x\right)  \right\vert =\left(  n+1\right)  .$ Now
using Lemma (\ref{numerical}) we have%
\begin{align*}
&  \left\vert f_{H}\left(  x|q\right)  -\frac{\sqrt{1-q}}{2\pi}\sqrt{4-\left(
1-q\right)  x^{2}}\left(  1+\sum_{k=2}^{n-1}(-1)^{k+1}q^{\binom{k}{2}}%
U_{2k-2}\left(  x\sqrt{1-q}/2\right)  \right)  \right\vert \\
&  =\left\vert \frac{\sqrt{1-q}}{2\pi}\sqrt{4-\left(  1-q\right)  x^{2}}%
\sum_{k\geq n}(-1)^{k+1}q^{\binom{k}{2}}U_{2k-2}\left(  x\sqrt{1-q}/2\right)
\right\vert \\
&  \leq\frac{\sqrt{1-q}}{2\pi}\sqrt{4-\left(  1-q\right)  x^{2}}\sum_{k\geq
n}\left\vert q\right\vert ^{\binom{k}{2}}\left(  2k-1\right)  \leq\frac
{\sqrt{1-q}}{\pi}\frac{2n\left\vert q\right\vert ^{n\left(  n-1\right)  /2}%
}{(1-q^{2})^{2}}.
\end{align*}

$ii)$ $\left\vert F_{H}\left(  y\right)  -\frac{1}{2}-\frac{1}{\pi}%
\arcsin\left(  \frac{y\sqrt{1-q}}{2}\right)  -\frac{\sqrt{4-\left(
1-q\right)  y^{2}}}{2\pi}\sum_{k=1}^{n-1}\left(  -1\right)  ^{k-1}q^{\binom
{k}{2}}U_{2n-2}\left(  \frac{y\sqrt{1-q}}{2}\right)  \right\vert
\allowbreak\leq\allowbreak\frac{\sqrt{4-\left(  1-q\right)  y^{2}}}{2\pi}%
\sum_{k=n}^{\infty}\left\vert q\right\vert ^{\binom{k}{2}}\left\vert
U_{2n-2}\left(  \frac{y\sqrt{1-q}}{2}\right)  \right\vert \allowbreak
\leq\allowbreak\frac{1}{2\pi}\sum_{k=n}^{\infty}\left\vert q\right\vert
^{\binom{k}{2}}\left(  1+\left\vert q\right\vert ^{k}\right)  $ \newline since
$\sup_{\left\vert x\right\vert \leq1}\sqrt{1-x^{2}}\left\vert U_{n}\left(
x\right)  \right\vert \allowbreak=\allowbreak1.$ Now to get $ii) $ we use
routine transformations and sum two geometric series.
\end{proof}

\begin{proof}
[Proof of Lemma \ref{symul}]Notice that comparing definitions of $f_{H}$ and
$f_{E}$ we have $f_{H}\left(  x|q\right)  \allowbreak=\allowbreak f_{E}\left(
x|q\right)  \allowbreak\times\allowbreak\frac{1-q^{9}}{1-q}\frac{1-q^{5}}%
{1-q}\allowbreak\times\allowbreak\prod_{k=4}^{\infty}\left(  \left(
1+q^{k}\right)  ^{2}-(1-q)x^{2}q^{k}\right)  \allowbreak\times\allowbreak
\prod_{k=1}^{\infty}\left(  1-q^{k}\right)  \allowbreak=\allowbreak
(1+q)(1-q^{3})\left(  1-q^{5}\right)  \left(  1-q^{9}\right)  \allowbreak
\times\allowbreak\prod_{k=4}^{\infty}\left(  \left(  1+q^{k}\right)
^{2}-(1-q)x^{2}q^{k}\right)  \left(  1-q^{k}\right)  .$ Now if $q\in
\lbrack0,1)$ we have
\begin{align*}
\left(  \left(  1+q^{k}\right)  ^{2}-(1-q)x^{2}q^{k}\right)  \left(
1-q^{k}\right)   &  \leq\left(  1+q^{k}\right)  ^{2}\left(  1-q^{k}\right) \\
&  =\left(  1-q^{2k}\right)  \left(  1+q^{k}\right)  .
\end{align*}
If $q\in(-1,0)$ then%
\begin{gather*}
\left(  \left(  1+q^{k}\right)  ^{2}-(1-q)x^{2}q^{k}\right)  \left(
1-q^{k}\right)  \leq\\
\left(  1-q^{k}\right)  \left\{
\begin{array}
[c]{ccc}%
\left(  \left(  1+q^{k}\right)  ^{2}-4q^{k}\right)  & if & k\text{ is odd}\\
\left(  1+q^{k}\right)  ^{2} & if & k\text{ is even}%
\end{array}
\right. \\
=\left(  1-q^{k}\right)  \left(  1+\left\vert q\right\vert ^{k}\right)  ^{2},
\end{gather*}
since then $\underset{\left\vert x\right\vert \leq2/\sqrt{1-q}}{\sup}\left(
\left(  1+q^{k}\right)  ^{2}-(1-q)x^{2}q^{k}\right)  \allowbreak
=\allowbreak\left(  \left(  1+q^{k}\right)  ^{2}-4q^{k}\right)  \allowbreak
=\allowbreak\left(  1-q^{k}\right)  ^{2}.$
\end{proof}

\subsection{Appendix}

Program in Mathematica that generates i.i.d. sequences from $f_{H}$.

\emph{QN[q\_,M\_]:=(Label[pocz];Y=y/.FindRoot[F[y,q]-RandomReal[],\{y,0\}];}

\emph{u=RandomReal[];t=newMM[q,M]/R[Y,q,M];If[t u%
$<$%
=1,Y,Goto[pocz]]);}

However it requires definition of function $F$ which is in fact function
$F_{E}$ of this paper. It is quite lengthy. newMM denotes function $M$ of this
paper. Further function $R$ denotes the ratio $f_{H}/f_{E}.$ Parameter $M$
denotes number that we insert instead of $\infty$ in the above mentioned
formulae. The above procedure produces $1$ observation from $f_{H}.$

Now \emph{AA[q\_,M\_,h\_]:=ListPlot[Table[QN[q,M],\{2000\}],PlotStyle-%
$>$%
Hue[h]];} produces table of $2000$ observation from $f_{H}(.|q) $ and plots it
in color $h.$ Then \emph{AA[.8,100,.4]} and \emph{AA[-.8,100,0]} produce plots
for $q\allowbreak=\allowbreak.8$ in color.$4$(green) and $q\allowbreak
=\allowbreak-.8$ in color $0$(red).

\begin{acknowledgement}
The author would like to thank all three referees for many sugesttions that
helped to improve the paper.
\end{acknowledgement}


\begin{thebibliography}{99}                                                                                               %


\bibitem[1]{AAR}Andrews, G.~E., Askey, R. and Roy, R.: 1999, \emph{Special
Functions}, Cambridge University Press, Cambridge, U.K.

\bibitem[2]{Ans04}Michael Anshelevich, ;2004, \emph{\ $q$-Levy processes}, J.
Reine Angew. Math. 576 (2004), 181-207.

\bibitem[3]{Bell}Bell, W., 1968, \emph{Special Functions}, D. van Nosrtand
Company Ltd, London.

\bibitem[4]{BKS97}Bo{\.{z}}ejko, M., K{\"{u}}mmerer, B., Speicher, R. : 1997,
\emph{$q$-Gaussian Processes: Non-commutative and Classical Aspects.}, Comm.
Math. Phys. \textbf{185(4(1))}, ~129--154.

\bibitem[5]{MBJW}Bo{\.{z}}ejko, M., Wysocza\'{n}ski, J. : 2001, \emph{Remarks
on }$t-$\emph{transformations of measures and convolutions}, Ann. Math. Instit
Poincar\'{e} Probab. Stat., \textbf{37}(6), 737-761.\textbf{6}(9), 4743-4756

\bibitem[6]{BB05}Bryc, W., Bo{\.{z}}ejko, M. : 2006, \emph{On a class of free
Levy laws related to a regression problem}, Journal of Functional Analysis
Volume 236 (2006), 59-77.

\bibitem[7]{BryMaWe}Bryc, W., Matysiak, W., Weso\l owski, J. , \emph{The bi -
Poisson process: a quadratic harness.} Annals of Probability 36 (2) (2008), s. 623-646

\bibitem[8]{bryc1}Bryc, W. (2001) Stationary random fields with linear
regressions. Annals of Probability 29, No. 1, 504-519.

\bibitem[8]{bryc2}Bryc, W. (2001) \emph{Stationary Markov chains with linear
regressions. Stochastic Processes and Applications} 93, 339-348.

\bibitem[9]{Bryc et al.}Bryc, W., Matysiak, W., Szab{\l }owski, P.~J. : 2005,
\emph{Probabilistic Aspects of al-Salam--Chihara Polynomials}, Proc. Amer.
Math. Soc. \textbf{\ 133},~1127--1134.

\bibitem[10]{Devroy}Devroye, L.: 1986, \emph{Non-Uniform Random Number
Generation}, Springer Verlag, New York.

\bibitem[11]{HLHM}van Leeuven, H., Maassen, H.: 1995, \emph{A }$\emph{q-}%
$\emph{Deformation of the Gauss distribution, }J. Math. Physics, \textbf{3}

\bibitem[12]{Szab}Szab\l owski, P.J. 2008 \emph{Probabilistic Implications of
symmetries of }$q$\emph{-Hermite and Al-Salam -Chihara Polynomials}, Infinite
Dimensional Analysis, Quantum Probability and Related Topic, 11(4), 513-522

\bibitem[13]{TMNT2008}W. Thistleton, W., Marsh, J. A., Nelson, K., Tsallis,C.,
2008, \emph{Generalized Box-Muller method for generating q-Gaussian random
deviates, }http://www.citeulike.org/user/orahcio/article/2859929

\bibitem[14]{Whittaker}Whittaker, E. T. and Watson, G. N. \emph{A course of
Modern Analysis, }Cambridge Univ Press, 1946 IV ed.
\end{thebibliography}
\end{document}